\documentclass [12pt, reqno]{amsart} 
\usepackage{times}
\usepackage{amsmath,amssymb,amsfonts,latexsym,eucal}
\usepackage[dvipsnames,svgnames,table]{xcolor}
\usepackage{enumitem}
\usepackage{tikz}

\allowdisplaybreaks[4] 

\newcommand\hlight[1]{\tikz[overlay, remember picture,baseline=-\the\dimexpr\fontdimen22\textfont2\relax]\node[rectangle,fill=white!50,rounded corners,fill opacity = 0.2,draw,thick,text opacity =1] {$#1$};}

\usepackage[bookmarksnumbered, colorlinks, linkcolor=red, citecolor=blue, plainpages]{hyperref}

\newtheoremstyle{myremark}     {10pt}{10pt}{}{}{\bfseries}{.}{.5em}{}

\newtheorem{thm}{Theorem}[section]

\newtheorem{lem}[thm]{Lemma}
\newtheorem{prop}[thm]{Proposition}
\theoremstyle{definition}

\theoremstyle{myremark}
\newtheorem{rem}[thm]{Remark}

\newtheorem{ques}[thm]{Question}

\newcommand{\abs}[1]{\left\vert#1\right\vert}

\newcommand{\norm}[1]{\left\Vert#1\right\Vert}

\DeclareMathOperator{\id}{id}
\let\ker\relax                                        
\DeclareMathOperator{\ker}{Ker}

\DeclareMathOperator{\Tr}{Tr}

\begin{document}

\title{A characterization of purely atomic von Neumann algebras}

\author[De]{Debabrata De}
\address{\hskip-\parindent
Debabrata De, Department of Mathematics, Indian Institute of Technology Guwahati, Guwahati 781039, India.}
	
\email{debabratade.ju@iitg.ac.in; debabratade.ju@gmail.com}

\subjclass[2020]{46L36, 46L51 $($primary$)$, 46L10, 46L53, 47B10 $($secondary$)$}

\keywords{von Neumann algebras, noncommutative Orlicz spaces}

\begin{abstract}
Let $M$ be a semifinite von Neumann algebra with a faithful normal semifinite trace $\tau$, and let $\Phi$
be an Orlicz function. We prove that the noncommutative Orlicz space $L^\Phi(M,\tau)$ is contained in $M$
if and only if $M$ is purely atomic and there is $\delta>0$ such that $\tau(e)\geq\delta$ for every nonzero
projection $e\in M$. We also prove that every left $M$-module map on
$L^\Phi(M,\tau)$ is bounded. As a consequence, a left ideal of $M$ contained in $L^\Phi(M,\tau)$ and dense
in it admits an $M$-invariant algebraic complement if and only if it is equal to $L^\Phi(M,\tau)$. These results
extend a theorem of Fidaleo and Zsid\'o from $L^2$ to all noncommutative Orlicz spaces.
\end{abstract}


\maketitle

\section{Introduction}\label{Sec:Intro} 
In the study of operator algebras, particularly in the theory of von Neumann algebras, factors serve as
building blocks from which all von Neumann algebras are constructed. There are several characterizations 
in the literature to determine the type of a von Neumann factor. Among all von Neumann factors, 
type $\mathrm{I}$ factors $($which are isomorphic to $\mathbf{B}(\mathcal{H})$ for some Hilbert space 
$\mathcal{H})$ have always been focal points in both operator algebras and operator theory. 
A von Neumann algebra is called \textit{purely atomic} if it is a direct sum of type $\mathrm{I}$ 
factors.  A natural question, inspired by the work of Buchholz, D'Antoni and Longo \cite{BDL90}, is the following:
\begin{ques}
Can purely atomic von Neumann algebras be characterized using their noncommutative $L^p$-spaces or generally by noncommutative Orlicz spaces, rather than only their GNS (that is, $L^2$) spaces?
\end{ques}
\noindent In this paper, we answer this question using noncommutative
$L^p$-spaces and, more generally, noncommutative Orlicz spaces,  associated with semifinite von Neumann algebras
with separable preduals, equipped with faithful normal
semifinite traces, for every $1\leq p<\infty$.

In the study of Kadison's ortho-unitary basis problem \cite{Ka67}, it is to be noted that 
GNS spaces $($with respect to faithful normal states$)$ of type $\mathrm{I}$ factors do not have both-sided uniformly bounded orthonormal bases \cite[Theorem 7.14]{DM23}. These limitations are due to the compactness of symmetric embeddings of von Neumann algebras into their GNS spaces \cite[Theorem 7.13]{DM23}. Recall that the symmetric embedding of a von Neumann algebra is a map that sends positive operators in the von Neumann algebra to elements in the canonical positive cone in the GNS space. Such an embedding was first investigated by D. Buchholz, C. D'Antoni and R. Longo, and they showed that the compactness of such an embedding characterizes type $\mathrm{I}$ factors \cite[Corollary 2.9]{BDL90}.  Recently, several equivalent criteria were established to determine the type of a factor, which involve the symmetric embedding of a von Neumann algebra into any of its bimodules induced by Markov maps \cite{DM21}, in terms of modular data and compactness of certain sets in bimodules induced by Markov maps on von Neumann algebras \cite{DM22} and the existence of both-sided uniformly bounded orthonormal bases in GNS spaces of von Neumann algebras \cite{DM23}.

Suppose $M\subseteq \mathbf{B}(\mathcal{H})$ is a $\sigma$-finite von Neumann factor represented on the GNS Hilbert space $\mathcal{H}$ $($with respect to a faithful normal state$)$ in standard form \cite{Ha75}. In a different direction, recently, K. H. Han, S. H. Kye and E. St\o rmer show that $M$ is of type $\mathrm{I}$ if and only if every normal completely bounded linear map $\phi$ from $M$ to $M$ has a certain type of normal extension $\tilde{\phi}$ on $\mathbf{B}(\mathcal{H})$ $($see \cite[Theorem 4.1]{HKS24}$)$.

All of the results above have two limitations. First, they use a faithful normal \emph{state}, so the 
algebra must be $\sigma$-finite. Second, and more importantly, all the aforesaid results are formulated in 
terms of $L^2$-spaces; that is, $L^2$-spaces are consistently present in all the above results.
Fidaleo and Zsid\'o removed the first limitation by working with faithful normal semifinite $(f.n.s.)$ 
weights. They proved an automatic continuity theorem for von Neumann algebras in standard form \cite{Ha75} 
and used it to characterize type $\mathrm{I}$ factors \cite[Proposition 6.2 and Corollary 6.3]{FZ16}. 
However, their characterization still relies on $L^2$-spaces. This leads to the following question: \emph{Is 
the Hilbert space structure of $L^2$ really needed for these characterizations?} For $p\neq 2$, 
the $L^p$-space is not a Hilbert space, so Hilbert space methods involving orthogonal projections are no 
longer available. So it is not clear a priori whether these characterizations extend to noncommutative 
$L^p$-spaces and noncommutative Orlicz spaces. Since we work with von Neumann algebras with separable predual, they are automatically 
$\sigma$-finite. We do not claim to remove the first limitation described
above. In this paper, we address the second limitation.

\subsection*{Main Result}	
In this paper, we extend a result of Fidaleo and Zsid\'o \cite{FZ16} to characterize purely atomic von 
Neumann algebras using noncommutative $L^p$-spaces for every $1\leq p<\infty$, and noncommutative Orlicz 
spaces $L^\Phi(M,\tau)$ (see Section \ref{Sec:MainResult}), showing that the Hilbert space structure is not 
essential. Throughout this paper, we consider noncommutative $L^p$-spaces and noncommutative Orlicz spaces 
associated with a semifinite von Neumann algebra $M$ with separable predual, equipped with a
faithful normal semifinite trace $\tau$. Let 
$\Phi:[0,\infty)\rightarrow [0,\infty)$
be a continuous convex function such that $\Phi(0)=0$ and $\Phi(t)>0$ for every $t>0$. Such a function $\Phi$ is called an \emph{Orlicz function}.
We write $\mathfrak{N}_\Phi=M\cap L^\Phi(M,\tau)$. We are now ready to state the main 
result of this paper.

\begin{thm}\label{Thm:PurelyAtomicVonNeumannAlgebras}
Let $M$ be a semifinite von Neumann algebra with separable predual, equipped with a $f.n.s.$ trace $\tau$ such that $\tau(1)=\infty$. Let $\Phi$ be an Orlicz function and let $\mathfrak{N}_\Phi:=M\cap L^\Phi(M,\tau)$. Let $\mathcal{X}\subseteq\mathfrak{N}_\Phi$ be a left ideal of $M$ which is dense in $L^\Phi(M,\tau)$ $($such an ideal exists, for example $\mathcal{X}=\mathfrak{N}_\Phi$, whenever $\Phi$ satisfies the $\Delta_2$-condition$)$. 
Then the following are equivalent:
\begin{enumerate}
\item $\mathcal{X}$ admits an $M$-invariant algebraic complement in $L^\Phi(M,\tau)$, that is, there exists a linear subspace $\mathcal{Y}\subseteq L^\Phi(M,\tau)$ such that 
\begin{align*}
L^\Phi(M,\tau)=\mathcal{X}+\mathcal{Y}, \quad \mathcal{X}\cap\mathcal{Y}=\{0\}, \quad\text{and}\quad a\mathcal{Y}\subseteq\mathcal{Y}\,\text{ for all } a\in M.
\end{align*}

\item $\mathcal{X}=\mathfrak{N}_\Phi$ and $\mathfrak{N}_\Phi=L^\Phi(M,\tau)$.
			
\item 	$\mathcal{X}=\mathfrak{N}_\Phi$ and $\underset{x\neq 0}{\underset{0\leq x\in M}{\sup}}\frac{\norm{x}_\infty}{\norm{x}_\Phi}<\infty$.
			
\item 	$\mathcal{X}=\mathfrak{N}_\Phi$  and the reduced von Neumann algebra $eMe$ is finite-dimensional for every projection $e\in\mathcal{P}(M)$ with $\tau(e)<\infty$.
			
\item 	$\mathcal{X}=\mathfrak{N}_\Phi$  and $M$ is purely atomic, that is, $M=\underset{i=1}{\overset{N}\oplus} M_i$, $N\leq \aleph_0$ with each $M_i$ a type $\mathrm{I}$ factor, and 
\begin{align*}
\tau(x)\geq \lambda (\underset{i=1}{\overset{N}\oplus} \Tr_i)(x), \quad 0\leq x\in M,
\end{align*}
holds true for some $\lambda>0$, where $\Tr_i$ is the canonical trace on $M_i$.
\end{enumerate}
\end{thm}

\subsection*{Methodology}
We now outline the methods used in this paper.\begin{enumerate}
\item[(a)] \emph{Automatic continuity.} We replace the automatic continuity result of \cite{FZ16} for $L^2$ by an automatic continuity theorem on noncommutative Orlicz spaces, in which we adapt the proof of \cite[Theorem 2.7]{AEGV21}: every left $M$-module map on $L^\Phi(M,\tau)$ is bounded $($Theorem \ref{Thm:AutoContOrlicz}$)$. In particular, this holds on $L^p(M,\tau)$ for every $1\leq p<\infty$.

\item[(b)] \emph{Closed graph theorem and the measure topology.} We use the structure theory of von Neumann algebras, the closed graph theorem, and the measure topology to prove the inclusion $L^\Phi(M,\tau)\subseteq M$ is equivalent $M$ is purely atomic and $\tau$ is bounded below on nonzero projections.

\item[(c)] \emph{Evaluating on projections.}
For every nonzero projection $e$ of finite trace,
$\norm{e}_\infty=1$ and $\norm{e}_p=\tau(e)^{1/p}$.
Thus, the inequality $\norm{x}_\infty\leq C\norm{x}_p$
implies that $\tau(e)\geq C^{-p}$ for every nonzero
projection $e$. Conversely, a uniform lower bound
$\delta>0$ on the traces of nonzero projections yields
$\norm{x}_\infty\leq\delta^{-1/p}\norm{x}_p$ by a
Chebyshev-type inequality. A similar argument applies
to noncommutative Orlicz spaces.
\end{enumerate}

\subsection*{Novelty of our work}
We now describe the main novelty of this paper.
\begin{enumerate}
\item[(a)] We extend the characterization of \cite[Proposition 6.2 and Corollary 6.3]{FZ16} from $L^2$ to all noncommutative Orlicz spaces. In particular, it is true for all noncommutative $L^p$-spaces, for $1\leq p<\infty$.

\item[(b)] We introduce the symmetric embedding $j_\Phi$ of a von Neumann algebra into noncommutative Orlicz spaces for an $f.n.s.$ trace, and show that the property that $j_\Phi$ is bounded below $($equivalently, onto$)$ characterizes purely atomic von Neumann algebras and, in the tracial state case, it characterizes finite-dimensional algebras (Theorem \ref{Thm:PurelyAtomic} and Remark \ref{Rem:IndepProof}$(4)$).

\item[(c)]  We give a new and self-contained proof of $(4)\Leftrightarrow (5)$ of Theorem \ref{Thm:PurelyAtomicVonNeumannAlgebras}, which is independent of the proof in \cite{FZ16} (see Remark \ref{Rem:IndepProof}). It also gives a noncommutative analog of the classical fact that $L^p(\mu)\subseteq L^\infty(\mu)$ if and only if $\mu$ is purely atomic with atoms of measure bounded away from zero \cite{Ro83, Vi85}.

\item[(d)] We show that a single, simple condition, namely that the traces of nonzero projections are bounded away from zero, controls everything: it is equivalent to $L^p(M,\tau)\subseteq M$ for one $($equivalently, all$)$ $1\leq p<\infty$, and to $L^\Phi(M,\tau)\subseteq M$ for one $($equivalently, all$)$ Orlicz functions $\Phi$. We also obtain the optimal constant in the inequality $\norm{x}_\infty\leq C\norm{x}_\Phi$ (Theorem \ref{Thm:PurelyAtomic} and Remark \ref{Rem:IndepProof}(5)).

\item[(e)] We prove that every left $M$-module map on a noncommutative Orlicz space $L^\Phi(M,\tau)$ is automatically bounded $($Theorem \ref{Thm:AutoContOrlicz}$)$. This extends \cite[Corollary 2.8]{AEGV21} beyond $L^p$-spaces, and this proof works for every symmetric Banach space of $\tau$-measurable operators.
\end{enumerate}

\subsection*{Organization of the paper}
The paper is organized as follows. In Section \ref{Sec:Preliminary}, we recall the basic facts on noncommutative $L^p$-spaces and noncommutative Orlicz spaces. In Section \ref{Sec:MainResult}, Theorem \ref{Thm:PurelyAtomic} characterizes purely atomic von Neumann algebras by the inclusion $L^p(M,\tau)\subseteq M$ and $L^\Phi(M,\tau)\subseteq M$, and the boundedness below of the symmetric embedding, in terms of the atomic structure of $M$ and the trace $\tau$.
In Theorem \ref{Thm:AutoContOrlicz}, we prove the automatic continuity theorem on $L^\Phi(M,\tau)$. We then prove Theorem \ref{Thm:PurelyAtomicVonNeumannAlgebras}, which gives the main characterization of purely atomic von Neumann algebras in this paper. In Theorem \ref{Thm:PurelyAtomicVonNeumannAlgebras3}, we deal with the tracial state case.

\section{Preliminaries}\label{Sec:Preliminary}
From now onwards, throughout this paper, unless otherwise mentioned, we denote $M\subseteq\mathbf{B}(\mathcal{H})$ as a semifinite von Neumann algebra with separable predual and $\tau$ as a faithful normal semifinite ($f.n.s.$ in short) trace on $M$. The norm on $M$ will be denoted by $\norm{\cdot}_\infty$. For $x\in L^0(M,\tau)$, we denote by $l(x)\in M$ the left support of $x$, that is, the projection onto $\overline{\text{range}(x)}$. Before we proceed, let us recall some basic facts on noncommutative $L^p$-spaces and noncommutative Orlicz spaces associated with the pair $(M,\tau)$. For more detailed theory on noncommutative $L^p$-spaces and noncommutative Orlicz spaces, we refer to \cite{Hi21, FK86, GL20}, and for the theory of von Neumann algebras, we refer to \cite{SZ19, KR97}.
\subsection{Noncommutative $L^p$-spaces}\label{SubSec:Preliminary}
 Let $L^0(M,\tau)$ denote the collection of $\tau$-measurable operators. 
 Then one can naturally extend $\tau$ from $M_+$, the positive part of $M$, to $L^0(M,\tau)_+$, the positive part of $L^0(M,\tau)$,  
 by setting: 
\begin{align*}
\tau(x)=\int_{0}^{\infty} \lambda\,\, d\tau(e_\lambda), \quad x\in L^0(M,\tau)_+,
\end{align*}
where $x=\int_{0}^{\infty} \lambda \,\,d e_\lambda$ denotes the spectral decomposition of $x$. For $1\leq p<\infty$ and $x\in L^0(M,\tau)$, we define $\norm{x}_{p}=\tau(\abs{x}^p)^{1/p}$. Then the noncommutative $L^p$-space $L^p(M,\tau)$ associated with the pair $(M,\tau)$ is defined as
\begin{align*}
L^p(M,\tau)=\{x\in L^0(M,\tau): \norm{x}_p<\infty\}.
\end{align*}
Note that $L^p(M,\tau)$ equipped with the norm $\norm{\cdot}_{p}$ becomes a Banach space, for $1\leq p<\infty$.

\noindent For  $x\in L^0(M,\tau)$ and $t>0$,  the $(t$-th$)$ generalized $s$-number $\mu_t(x)$ is defined by 
\begin{align*}
\mu_t(x)=\inf\{s>0:\tau(e_{(s,\infty)}(\abs{x}))\leq t\},
\end{align*}
where $e_{(s,\infty)}(\abs{x})$ is the spectral projection of $\abs{x}$. It is known that for an element $x\in L^0(M,\tau)$, we have 
\begin{align*}
\tau(\abs{x}^p)=\int_{0}^{\infty}\mu_t(x)^p \, dt.
\end{align*}
Note that $M\cap L^p(M,\tau)$ is dense in $ L^p(M,\tau)$ with respect to  $\norm{\cdot}_{p}$ for $1\leq p<\infty$.
	
\subsection{Noncommutative Orlicz spaces}\label{Subsec:Orlicz}
Recall that an \emph{Orlicz function} is a
continuous convex function $\Phi:[0,\infty)\rightarrow [0,\infty)$
such that $\Phi(0)=0$ and $\Phi(t)>0$ for every $t>0$.
Note that $\Phi$ is strictly increasing
and $\Phi(t)\rightarrow\infty$ as $t\rightarrow\infty$. Thus, $\Phi$ has a continuous, strictly increasing inverse
$\Phi^{-1}:[0,\infty)\rightarrow [0,\infty)$. The noncommutative Orlicz space $L^\Phi(M,\tau)$ associated with $(M,\tau)$ and an Orlicz function $\Phi$ is defined as the collection of all
$x\in L^0(M, \tau)$ for which there exists $\alpha>0$ such that 
\begin{align*}
\tau(\Phi(\alpha\abs{x}))<\infty.    
\end{align*}
Then $L^\Phi(M,\tau)$ equipped with the Luxemburg-Nakano norm:
\begin{align*}
\norm{x}_\Phi=\inf\left\{\alpha>0: \tau\left(\Phi\left(\frac{\abs{x}}{\alpha}\right)\right)\leq 1\right\}, \quad \quad x\in L^\Phi(M, \tau),    
\end{align*}
becomes a Banach space \cite[Theorem 5.45]{GL20}.
In particular, for $\Phi(t)=t^p$, $(L^\Phi(M,\tau),\norm{\cdot}_\Phi)$ coincides with $(L^p(M,\tau), \norm{\cdot}_p)$.

We begin with the following elementary lemma.
\begin{lem}\label{Lem:Orlicz}
Let $\Phi$ be an Orlicz function.
\begin{enumerate}
\item If $0\neq x\in L^\Phi(M,\tau)$ and $s>0$, then 
\begin{align*}
\tau\left(e_{(s,\infty)}(\abs{x})\right)\leq \frac{1}{\Phi\left(s/\norm{x}_\Phi\right)}.
\end{align*}
In particular, convergence in $\norm{\cdot}_\Phi$ implies convergence in the measure topology.
\item If $e\in\mathcal{P}(M)$ is nonzero with $\tau(e)<\infty$, then $e\in L^\Phi(M,\tau)$ and $\norm{e}_\Phi=1/\Phi^{-1}\left(1/\tau(e)\right)$.
\end{enumerate}
\end{lem}

\begin{proof}
\noindent $(1)$. Put $e_s=e_{(s,\infty)}(\abs{x})$. Since $\Phi$ is increasing, $\Phi(\frac{t}{\lambda})\geq\Phi(\frac{s}{\lambda})\chi_{(s,\infty)}(t)$ for $t\geq 0$. By functional calculus, we have $\Phi(\frac{\abs{x}}{\lambda})\geq \Phi(\frac{s}{\lambda})\,e_s$. Note that if $\lambda>\norm{x}_\Phi$, then 
$\tau(\Phi(\frac{\abs{x}}{\lambda}))\leq 1$. Hence, $\Phi(\frac{s}{\lambda})\,\tau(e_s)\leq 1$ for all $\lambda>\norm{x}_\Phi$. By letting $\lambda\downarrow\norm{x}_\Phi$ and using the continuity of $\Phi$, we get the claim. 

If $\norm{x_n-x}_\Phi\rightarrow 0$, then for each $s>0$, we have 
\begin{align*}
\tau(e_{(s,\infty)}(\abs{x_n-x}))\leq \frac{1}{\Phi(\frac{s}{\norm{x_n-x}_\Phi})}\rightarrow 0,    
\end{align*}
since $\Phi(t)\rightarrow\infty$ as $t\rightarrow\infty$. This proves the convergence in measure.

\noindent $(2)$. Since $\Phi(0)=0$, by functional calculus, $\Phi(\frac{e}{\lambda})=\Phi(\frac{1}{\lambda})\,e$. Thus, $\tau(\Phi(\frac{e}{\lambda}))=\Phi(\frac{1}{\lambda})\tau(e)$. Hence, 
\begin{align*}
\tau(\Phi(\frac{e}{\lambda}))=\Phi(\frac{1}{\lambda})\tau(e)\leq 1 \quad\text{ if and only if }  \quad \frac{1}{\lambda}\leq\Phi^{-1}\left(\frac{1}{\tau(e)}\right), 
\end{align*}
which gives the formula for $\norm{e}_\Phi$. This completes the proof.
\end{proof}

\section{Main Result}\label{Sec:MainResult}
This is the main section of this paper. In this section, we characterize purely atomic von Neumann algebras among semifinite ones with respect to symmetric embeddings of von Neumann algebras into noncommutative Orlicz spaces.
Our result also extends a result of \cite{FZ16}  in noncommutative $L^p$-spaces from the case of $p=2$ to all $1\leq p<\infty$.

We also use the symmetric embedding of $M$ into $L^\Phi(M,\tau)$. For an Orlicz function $\Phi$, put $\mathfrak{N}_\Phi:=M\cap L^\Phi(M,\tau)$ and consider the map
\begin{align*}
j_\Phi: (\mathfrak{N}_\Phi,\norm{\cdot}_\infty)\rightarrow L^\Phi(M,\tau),\qquad j_\Phi(x)=x.
\end{align*}
It is injective, has dense range, and maps the positive part of $\mathfrak{N}_\Phi$ into the positive cone of $L^\Phi(M,\tau)$. In particular, when $\Phi_p(t)=t^p$, we denote 
$\mathfrak{N}_{\Phi_p}$ and $j_{\Phi_p}$ as $\mathfrak{N}_p$ and $j_p$, respectively.

When $\tau$ is a faithful normal tracial state and $\Phi(t)=t^2$, $j_\Phi$ is the symmetric embedding $M\ni x\mapsto x\xi_\tau\in L^2(M,\tau)$ into the GNS space (see, \cite{BDL90}). We therefore call $j_\Phi$ the \emph{symmetric embedding} of $M$ into $L^\Phi(M,\tau)$. In \cite{BDL90}, characterization of type $\mathrm{I}$ factors is obtained by \emph{compactness} of the symmetric embedding. Here, instead of \emph{compactness}, we characterize it by the boundedness of its \emph{inverse}, that is, by $j_\Phi$ being bounded below. We begin with the characterization of the inclusion $L^\Phi(M,\tau)\subseteq M$ in terms of a purely atomic von Neumann algebra.
	
\begin{thm}\label{Thm:PurelyAtomic}
Let $M$ be a semifinite von Neumann algebra with separable predual, equipped with a $f.n.s.$ trace $\tau$. For $1\leq p<\infty$ and an Orlicz function $\Phi$, let $\mathfrak{N}_p$, $\mathfrak{N}_\Phi$ and $j_\Phi$ be as above.	Then the following statements are equivalent:
		
\begin{enumerate}
\item $\mathfrak{N}_p=L^p(M,\tau)$, that is, $L^p(M,\tau)\subseteq M$.
			
\item $\underset{x\neq 0}{\underset{0\leq x\in M}{\sup}}\frac{\norm{x}_\infty}{\norm{x}_p}<\infty$.
			
\item The reduced von Neumann algebra $eMe$ is finite-dimensional for every projection $e\in\mathcal{P}(M)$ with $\tau(e)<\infty$.
			
\item $M$ is purely atomic, that is, $M=\underset{i=1}{\overset{N}\oplus} M_i$, $N\leq \aleph_0$ with each $M_i$ a type $\mathrm{I}$ factor, and 
\begin{align}\label{Eq:AtomicTraceBound}
\tau(x)\geq \lambda (\underset{i=1}{\overset{N}\oplus} \Tr_i)(x), \quad 0\leq x\in M,
\end{align}
holds true for some $\lambda>0$, where $\Tr_i$ is the canonical trace on $M_i$.

\item There exists $\delta>0$ such that $\tau(e)\geq\delta$ for every nonzero projection $e\in\mathcal{P}(M)$.

\item $L^\Phi(M,\tau)\subseteq M$.

\item There exists $C>0$ such that $\norm{x}_\infty\leq C\norm{x}_\Phi$ for all $x\in M\cap L^\Phi(M,\tau)$.

\item The symmetric embedding $j_\Phi$ is bounded below, that is, there is $C>0$ with $\norm{x}_\infty\leq C\norm{j_\Phi(x)}_\Phi$ for all $x\in\mathfrak{N}_\Phi$.

\item The symmetric embedding $j_\Phi$ is surjective.
\end{enumerate}
Moreover, if condition $(5)$ holds, then 
\begin{align}\label{Eq:BestConstant}
\norm{x}_\infty\leq \Phi^{-1}(1/\delta)\,\norm{x}_\Phi,\quad x\in L^\Phi(M,\tau).
\end{align}
Consequently, conditions $(1)$, $(2)$, $(6)$, $(7)$, $(8)$ and $(9)$ hold for one $1\leq p<\infty$ and one Orlicz function $\Phi$ if and only if they hold for all of them.
\end{thm}

\begin{proof}
Recall that for Orlicz function $\Phi_p(t)=t^p$, $L^{\Phi_p}(M,\tau)=L^p(M,\tau)$, 
$\norm{\cdot}_{\Phi_p}=\norm{\cdot}_p$ and $\Phi_p^{-1}(1/\delta)=\delta^{-1/p}$ for $1\leq p<\infty$. 

Thus, condition $(1)$ is condition $(6)$ for $\Phi=\Phi_p$. Since $\norm{x}_\infty=\norm{\abs{x}}_\infty$, $\norm{x}_\Phi=\norm{\abs{x}}_\Phi$ and $\abs{x}\in\mathfrak{N}_\Phi$, for $x\in\mathfrak{N}_\Phi$, and since $\norm{x}_p=\infty$ for $0\leq x\in M\setminus L^p(M,\tau)$, condition $(2)$ is condition $(7)$ for $\Phi=\Phi_p$. Moreover, since $j_\Phi(x)=x$, condition $(8)$ is exactly condition $(7)$, and condition $(9)$ says that $\mathfrak{N}_\Phi=L^\Phi(M,\tau)$, which is condition $(6)$.

Note that conditions $(3)$, $(4)$ and $(5)$ do not involve $\Phi$.
Thus, to prove conditions $(1)$-$(9)$ are equivalent, it is sufficient to show that 
\begin{align*}
(5)\Rightarrow(4)\Rightarrow(3)\Rightarrow(5)\qquad\text{ and }\qquad(6)\Rightarrow(7)\Rightarrow(5)\Rightarrow(6) \text{ and } (5)\Rightarrow \text{ Eq. }\eqref{Eq:BestConstant}.
\end{align*}

\vspace{2mm}

\noindent $(6)\Rightarrow(7)$. By $(6)$, the inclusion map $\iota_\Phi\colon L^\Phi(M,\tau)\rightarrow M$ is well defined. Let $x_n\rightarrow x$ in $L^\Phi(M,\tau)$ and $x_n\rightarrow y$ in $M$. Then $x_n\rightarrow x$ in the measure topology by Lemma \ref{Lem:Orlicz}$(1)$, and $x_n\rightarrow y$ in the measure topology by \cite[Lemma 4.18(4)]{Hi21}. Since the measure topology is Hausdorff \cite[Theorem 4.12]{Hi21}, $x=y$. As $L^\Phi(M,\tau)$ and $M$ are Banach spaces, by the closed graph theorem, $\iota_\Phi$ is bounded, which proves $(7)$.

\vspace{2mm}

\noindent $(7)\Rightarrow(5)$. Let $0\neq e\in\mathcal{P}(M)$. If $\tau(e)=\infty$, there is nothing to prove. Thus, we can assume that $\tau(e)<\infty$. By Lemma \ref{Lem:Orlicz}$(2)$ and $(7)$,
\begin{align*}
1=\norm{e}_\infty\leq C\norm{e}_\Phi=\frac{C}{\Phi^{-1}(1/\tau(e))},
\end{align*}
so that $1/\tau(e)\leq\Phi(C)$. Hence, $(5)$ holds with $\delta=1/\Phi(C)$. 

\vspace{2mm}

\noindent $(5)\Rightarrow(6)$. Let $0\neq x\in L^\Phi(M,\tau)$ and $s>\Phi^{-1}(1/\delta)\norm{x}_\Phi$. Then $\Phi(s/\norm{x}_\Phi)>1/\delta$, so Lemma \ref{Lem:Orlicz}$(1)$ gives $\tau(e_{(s,\infty)}(\abs{x}))<\delta$. By $(5)$, $e_{(s,\infty)}(\abs{x})=0$, that is, $\abs{x}\leq s1$. Hence $\abs{x}$ is bounded, and since $x$ is affiliated with $M$, $x\in M$. Moreover, since $\norm{x}_\infty\leq s$, by letting $s\downarrow\Phi^{-1}(1/\delta)\norm{x}_\Phi$ gives Eq. \eqref{Eq:BestConstant}.

\vspace{2mm}

\noindent $(5)\Rightarrow(4)$. We prove this implication in three steps.

\noindent \emph{Step 1:} Let $0\neq e\in\mathcal{P}(M)$ with $\tau(e)<\infty$. In this step, we show that $eMe$ is finite-dimensional. 

\vspace{2mm}

Let $0\neq e\in\mathcal{P}(M)$ with $\tau(e)<\infty$, and suppose
that there exists $\delta>0$ such that $\tau(f)\geq\delta$ for every
nonzero projection $f\leq e$. Put $\mathcal{R}=eMe$. If $f_1,\ldots,f_n$ are
mutually orthogonal nonzero projections in $\mathcal{R}$, then
\begin{align*}
n\delta
\leq \sum_{k=1}^{n}\tau(f_k)
= \tau\left(\sum_{k=1}^{n}f_k\right)
\leq \tau(e).
\end{align*}
Thus every such family has at most
$m=\lfloor\tau(e)/\delta\rfloor$ members.

Choose a family $z_1, \ldots, z_r$ of mutually orthogonal nonzero
projections in $\mathcal{Z}(\mathcal{R})$ of largest possible cardinality. Then
$\sum_{j=1}^{r}z_j=e$, since otherwise the nonzero central projection
$e-\sum_{j=1}^{r}z_j$ could be added to the family. Moreover, each $z_j$ is minimal in $\mathcal{Z}(\mathcal{R})$. Indeed, if
$q\in\mathcal{P}(\mathcal{Z}(\mathcal{R}))$ satisfies $0<q<z_j$, then
replacing $z_j$ by $q$ and $z_j-q$ gives a family of $r+1$
mutually orthogonal nonzero central projections, contradicting the
choice of $r$. Hence,
$\mathcal{Z}(z_j\mathcal{R})=z_j\mathcal{Z}(\mathcal{R})=\mathbb{C}z_j$, so $z_j\mathcal{R}$
is a factor.

For each $j$, choose a family of mutually orthogonal nonzero projections in $z_j\mathcal{R}$ of largest possible cardinality. Each member
of this family is minimal in $z_j\mathcal{R}$, since a proper nonzero subprojection $q<p$ would allow us to replace $p$ by $q$ and $p-q$.
Thus, the factor $z_j\mathcal{R}$ contains a minimal projection. Consequently, $z_j\mathcal{R}$ is a type $\mathrm{I}$ factor. Thus,  
\begin{align*}
eMe=\bigoplus_{j=1}^{r}z_j\mathcal{R}.
\end{align*}
Since $\tau$ is faithful and $\tau(e)<\infty$, we have $\tau(z_j)<\infty$, so $z_j\mathcal{R}$ is finite dimensional. Hence, $eMe$ is finite dimensional.

\vspace{2mm}

\noindent\emph{Step 2:} In this step, we show that $M$ is atomic.

\vspace{2mm}

Let $0\neq g\in\mathcal{P}(M)$. Since $\tau$ is semifinite, there is a nonzero projection $e\leq g$ with $\tau(e)<\infty$. By Step 1, $eMe$ is finite-dimensional, so it contains a projection $q\neq0$ which is minimal in $eMe$. Note that $qMq=q(eMe)q=\mathbb{C}q$. Thus,  $q\leq g$ is a minimal projection of $M$. Hence, $M$ is atomic, and therefore, $M=\underset{i=1}{\overset{N}\oplus}M_i$ with each $M_i$ a type $\mathrm{I}$ factor \cite{SZ19}; moreover $N\leq\aleph_0$, since $M$ has separable predual.

\vspace{2mm}

\noindent\emph{Step 3:} In this step, we show the existence of $\lambda>0$ such that Eq. \eqref{Eq:AtomicTraceBound} holds.

\vspace{2mm} 
Fix $1\leq i\leq N$. Let $p_i$ be the central projection in $M$ with $Mp_i=M_i$. Note that $\tau_i=\tau(p_i\,\cdot)$ is a $f.n.s.$ trace on $M_i$. Then, there exists $\lambda_i>0$ such that $\tau_i=\lambda_i\Tr_i$ \cite{SZ19}.
If $q$ is a minimal projection of $M_i$, then $\lambda_i=\tau(q)$. By $(5)$, there exists $\delta>0$ such that $\tau(q)\geq\delta$. Hence, $\lambda_i\geq \delta$. Consequently,
\begin{align*}
\tau(x)=\sum_{i=1}^N\lambda_i\Tr_i(xp_i)\geq\delta\,(\underset{i=1}{\overset{N}\oplus}\Tr_i)(x),\qquad 0\leq x\in M,
\end{align*}
which is $(4)$ with $\lambda=\delta$.

\smallskip
\noindent $(4)\Rightarrow(3)$. Assume that $(4)$ holds. Let $e\in\mathcal{P}(M)$ with $\tau(e)<\infty$. There exist central projections $p_i$ in $M$ such that $ Mp_i = M_i$. Note that
\begin{align*}
\infty>\tau(e)\geq\lambda\sum_i\Tr_i(ep_i).
\end{align*}
Since $\Tr_i(ep_i)\geq 1$ whenever $ep_i\neq 0$, the set $F=\{i: e p_i\neq 0\}$ is finite and each $ep_i$ has finite rank. Consequently, 
$ep_i M_i p_i e$ is finite-dimensional.
Hence, $eMe=\oplus_{i\in F}\,ep_i M_i p_i e$ is finite-dimensional.

\vspace{2mm}

\noindent $(3)\Rightarrow(5)$. Suppose $(5)$ fails. Then there exist nonzero projections $f_k\in\mathcal{P}(M)$
such that $\tau(f_k)<2^{-k}$ for every $k\geq 1$.
Set $e=\bigvee_{k\geq 1}f_k$. Note that Kaplansky's formula,
\begin{align*}
(p\vee q)-q\sim p-(p\wedge q),
\end{align*}
implies that $\tau(p\vee q)\leq\tau(p)+\tau(q)$ for all projections $p,q\in\mathcal{P}(M)$. Hence, by induction and normality of $\tau$, we have
\begin{align*}
\tau(e)
&=\lim_{n\to\infty}\tau\left(\bigvee_{k=1}^n f_k\right)
\leq\sum_{k=1}^{\infty}\tau(f_k)
\leq 1.
\end{align*}
By $(3)$, we have
\begin{align*}
eMe\cong\bigoplus_{j=1}^m M_{n_j}(\mathbb{C}).
\end{align*}
Since the restriction of $\tau$ to $eMe$ is a faithful normal finite trace,
so it has the form
\begin{align*}
\tau|_{eMe}=\bigoplus_{j=1}^m c_j\Tr_{n_j},
\qquad c_j>0,
\end{align*}
where $\Tr_{n_j}$ denotes standard trace on $M_{n_j}(\mathbb{C})$ that takes value $1$ on each minimal projection.
Thus, every nonzero projection in $eMe$ has trace at least
$\min_{1\leq j\leq m}c_j>0$. This contradicts the fact that
$0\neq f_k\leq e$ and $\tau(f_k)\to0$.
This proves $(5)$.
\end{proof}

\begin{rem}\label{Rem:IndepProof}
\noindent $(1)$. The conditions $(3)$ and $(4)$ do not involve any noncommutative Orlicz spaces or $L^p$-spaces. Therefore, to establish this equivalence, one can refer to the proof in \cite[Proposition 6.2]{FZ16}. However, in this paper, to show $(3)\Leftrightarrow (4)$, we show that $(3)\Rightarrow (5)\Rightarrow (4)\Rightarrow (3)$, with the proof relying on the structure of von Neumann algebras, which is independent of the proof in \cite[Proposition 6.2]{FZ16}.

\vspace{2mm}

\noindent $(2)$. Note that Theorem \ref{Thm:PurelyAtomic} is a noncommutative analog of the classical result: 
For a measure space $(\Omega,\mu)$, $L^p(\mu)\subseteq L^\infty(\mu)$ if and only if $\mu$ is purely atomic and the measures of its atoms are bounded away from zero; see, e.g., \cite{Ro83, Vi85}.    

\noindent $(3)$. In this paper, we have not considered Haagerup
$L^p$-spaces $($with respect to a $f.n.s.$ weight$)$ because for Haagerup $L^p$-spaces $L^p(M)$, $L^{p_1}(M)\cap L^{p_2}(M)=\{0\}$ if $p_1\neq p_2$ with $p_1,p_2\in [1,\infty]$.

\noindent $(4)$. The assumption $\tau(1)=\infty$ is not used in the proof of
Theorem \ref{Thm:PurelyAtomic}. In particular, the proof also applies when
$\tau$ is a faithful normal tracial state. Then condition $(3)$ with $e=1$ shows that the conditions of Theorem \ref{Thm:PurelyAtomic} are equivalent to $M$ being finite-dimensional.
 
\noindent $(5)$. Let $\delta=\inf\{\tau(e):0\neq e\in\mathcal{P}(M)\}$. If $(4)$ holds and $\tau=\lambda_i\Tr_i$ on $M_i$, then $\tau(q)=\lambda_i$ for every minimal projection $q$ of $M_i$, while every nonzero projection $e$ satisfies $\tau(e)\geq\lambda_i\Tr_i(ep_i)\geq\lambda_i$ for any $i$ with $ep_i\neq 0$. Hence, $\delta=\inf_i\lambda_i$. By applying Eq. \eqref{Eq:BestConstant} with $\Phi(t)=t^p$ to minimal projections shows that the constant $\delta^{-1/p}$ in $\norm{x}_\infty\leq\delta^{-1/p}\norm{x}_p$ is optimal.
\end{rem}

We now establish an automatic continuity result for
$L^\Phi(M,\tau)$ analogous to \cite[Corollary 3.2]{FZ16}, in which we extend the result from $L^2$-spaces to noncommutative Orlicz spaces. Our proof is based on the proof of \cite[Theorem 2.7]{AEGV21}.

\begin{thm}\label{Thm:AutoContOrlicz}
Let $M$ be a semifinite von Neumann algebra equipped with an $f.n.s.$ trace $\tau$, and let $\Phi$ be an Orlicz function. Then every linear map $T: L^\Phi(M,\tau)\rightarrow L^\Phi(M,\tau)$ satisfying
\begin{align}\label{Eq:LeftMulti}
T(ax)=aT(x),\qquad a\in M,\ x\in L^\Phi(M,\tau),
\end{align}
is bounded. In particular, for $\Phi(t)=t^p$, the same holds on $L^p(M,\tau)$ for every $1\leq p<\infty$.
\end{thm}
 
\begin{proof}
To prove this result, we follow the proof of  \cite[Theorem 2.7]{AEGV21}.

Denote the separating space of $T$ by
\begin{align*}
\mathcal{S}(T)
=\{y\in L^\Phi(M,\tau): \text{there is }(x_n)\subseteq L^\Phi(M,\tau) \text{ with } x_n\rightarrow 0 \text{ and } T(x_n)\rightarrow y\}.    
\end{align*}
Note that $\mathcal{S}(T)$ is a closed subspace of $L^\Phi(M,\tau)$ \cite[Proposition 5.1.2]{Da00}. By the closed graph theorem, $T$ is bounded if and only if $\mathcal{S}(T)=\{0\}$. Thus, we prove $\mathcal{S}(T)=\{0\}$ in four steps. 

\vspace{2mm}

\noindent \emph{Step 1:} In this step, we show that $\mathcal{S}(T)$ is a left $M$-submodule of $L^\Phi(M,\tau)$. Indeed, let $a\in M$ and $y\in \mathcal{S}(T)$. There exists $(x_n)\subseteq L^\Phi(M,\tau)$ with $x_n\rightarrow 0$ and $T(x_n)\rightarrow y$.   
Since $ax_n\rightarrow 0$ and $T(ax_n)=aT(x_n)\rightarrow ay$, we have $M\mathcal{S}(T)\subseteq\mathcal{S}(T)$.

\noindent \emph{Step 2:} In this step, we show that there exists a central projection $z\in\mathcal{P}(\mathcal{Z}(M))$ such that
\begin{align*}
\mathcal{I}(T)=\{a\in M : a\mathcal{S}(T)=\{0\}\}=Mz.   
\end{align*}

Indeed, note that 
\begin{align*}
\mathcal{I}(T)=\{a\in M : a\mathcal{S}(T)=\{0\}\}   
\end{align*}
is a closed left ideal of $M$. By Step 1, it is also a right ideal. Indeed, if $a\in\mathcal{I}(T)$ and $b\in M$, then $ab \mathcal{S}(T)\subseteq a\mathcal{S}(T)=\{0\}$. 

We now show that $\mathcal{I}(T)$ is weak$^*$-closed in $M$. Let $y\in L^\Phi(M,\tau)$ and $a\in M$. Suppose $ay=0$. Then, $ay\xi=0$ for all $\xi\in\mathfrak{D}(y)$. Since $l(y)$ is the projection onto $\overline{\text{range}(y)}$, it follows that $a l(y)=0$. Conversely, if $a l(y)=0$, then $ay=al(y)y=0$. Thus, we have
\begin{align*}
\{a\in M: ay=0\}
=\{a\in M: a l(y)=0\}.
\end{align*}
Since $l(y)\in M$, it follows that $\{a\in M: a l(y)=0\}$ is weak$^*$-closed in $M$. Consequently, 
\begin{align*}
\mathcal{I}(T)
=\bigcap_{y\in\mathcal{S}(T)}\{a\in M: a\,l(y)=0\}  
\end{align*}
is weak$^*$-closed in $M$. By \cite[Theorem 6.8.8]{KR97}, there exists a central projection $z\in\mathcal{P}(\mathcal{Z}(M))$ such that
$\mathcal{I}(T)=Mz$.

\noindent\emph{Step 3:} In this step, we show that $(1-z)M$ is finite-dimensional. 

Suppose that $(1-z)M$ is infinite-dimensional. Then 
there exists a sequence of non-zero mutually orthogonal projections $(e_k)$ in $(1-z)M$ such that $e_k\leq 1-z$. 
For each $n\in\mathbb{N}$, set
\begin{align*}
f_n=\bigvee_{k\geq n}e_k.    
\end{align*}
Let $R_n: L^\Phi(M,\tau)\rightarrow L^\Phi(M,\tau)$ be the bounded operator defined by $R_n(x)=f_n x$, for $x\in L^\Phi(M,\tau)$. Since $T$ is a left module map, $TR_n-R_nT=0$. By the stability lemma \cite[Corollary 5.2.7]{Da00}, the decreasing sequence of closed subspaces
\begin{align*}
\overline{R_1 R_2\cdots R_n\,\mathcal{S}(T)}=\overline{f_n\,\mathcal{S}(T)}    
\end{align*}
stabilizes, that is, there is $N$ with $\overline{f_N\mathcal{S}(T)}=\overline{f_{N+1}\mathcal{S}(T)}$. Since $e_N f_N=e_N$, $e_N f_{N+1}=0$ and left multiplication by $e_N$ is continuous on $L^\Phi(M,\tau)$, we have
\begin{align*}
e_N\mathcal{S}(T)=e_Nf_N\mathcal{S}(T)\subseteq e_N\overline{f_{N+1}\mathcal{S}(T)}\subseteq\overline{e_Nf_{N+1}\mathcal{S}(T)}=\{0\}.    
\end{align*}
Thus, $0\neq e_N\in\mathcal{I}(T)=Mz$, which is contradiction as $e_N\leq 1-z$. Hence, $(1-z)M$ is finite-dimensional.
 
\vspace{2mm}
\noindent\emph{Step 4}: In this step, we show that 
$\mathcal{S}(T)=\{0\}$.

First, note that by Step 3, $(1-z) L^\Phi(M,\tau)$ is finite-dimensional. Thus, the linear map $x\mapsto T((1-z)x)=(1-z)T(x)$ is continuous on $L^\Phi(M,\tau)$. Let $y\in\mathcal{S}(T)$. Thus, there exist $(x_n)\subseteq L^\Phi(M,\tau)$ such that
$x_n\rightarrow 0$ and $T(x_n)\rightarrow y\in\mathcal{S}(T)$. Note that $(1-z)T(x_n)\rightarrow 0$ and $(1-z)T(x_n)\rightarrow(1-z)y$. Hence, $(1-z)\mathcal{S}(T)=\{0\}$. So $1-z\in\mathcal{I}(T)=Mz$, which forces $1-z=0$. Therefore $1\in\mathcal{I}(T)$, that is, $\mathcal{S}(T)=\{0\}$. Hence, $T$ is bounded. This completes the proof.
\end{proof}

\begin{rem}
The proof of Theorem \ref{Thm:AutoContOrlicz} and \cite[Theorem 2.7]{AEGV21} use only three facts: $(i)$ $\norm{axb}_\Phi\leq\norm{a}_\infty\norm{x}_\Phi\norm{b}_\infty$, $a,b\in M$ and $x\in L^\Phi(M,\tau)$; $(ii)$ $ay=0$ if and only if $a\,l(y)=0$; $(iii)$ $zL^\Phi(M,\tau)$ is finite-dimensional whenever $z$ is a central projection with $\dim zM<\infty$.
It therefore holds verbatim for every symmetric Banach space $E(M,\tau)$ of $\tau$-measurable operators.
\end{rem}

We are now ready to prove the main result (Theorem \ref{Thm:PurelyAtomicVonNeumannAlgebras}) of the paper, which characterizes purely atomic von Neumann algebras among semifinite von Neumann algebras.

\begin{proof}[Proof of Theorem \ref{Thm:PurelyAtomicVonNeumannAlgebras}]
$(1)\Rightarrow (2)$. Let $P: L^\Phi(M,\tau)\rightarrow L^\Phi(M,\tau)$ be the projection onto $\mathcal{X}$ along $\mathcal{Y}$, that is, $P(x+y)=x$ for $x\in\mathcal{X}$ and $y\in\mathcal{Y}$. Then, $P$ is a well-defined linear idempotent with range $\mathcal{X}$. Let $a\in M$, $x\in\mathcal{X}$ and $y\in\mathcal{Y}$. Since $\mathcal{X}$ is a left ideal of $M$, we have $ax\in\mathcal{X}$, and $ay\in\mathcal{Y}$ by the invariance of $\mathcal{Y}$. Hence, we have 
\begin{align*}
P\left(a(x+y)\right)=P(ax+ay)=ax=aP(x+y).
\end{align*}
Thus $P(az)=aP(z)$ for all $a\in M$ and $z\in L^\Phi(M,\tau)$. Consequently, $P$ is bounded by Theorem \ref{Thm:AutoContOrlicz}. Hence, $\mathcal{X}=\ker(\id-P)$ is closed in $L^\Phi(M,\tau)$. Since $\mathcal{X}$ is dense in $L^\Phi(M,\tau)$, we have $\mathcal{X}=L^\Phi(M,\tau)$. Since $\mathcal{X}\subseteq\mathfrak{N}_\Phi\subseteq L^\Phi(M,\tau)=\mathcal{X}$, we get $\mathcal{X}=\mathfrak{N}_\Phi=L^\Phi(M,\tau)$. Hence, $(2)$ holds.

\vspace{2mm}
\noindent $(2)\Rightarrow (1)$. The proof of $(2)\Rightarrow (1)$ is trivial. 

Now the equivalences of $(2)\Leftrightarrow (3) \Leftrightarrow (4)\Leftrightarrow (5)$ follow from Theorem \ref{Thm:PurelyAtomic}. This completes the proof.
\end{proof}

\begin{prop}\label{Prop:LeftIdealBound}
Let $M$, $\tau$, $\Phi$ and  $\mathcal{X}$ be as in Theorem \ref{Thm:PurelyAtomicVonNeumannAlgebras}. Then the following statements are equivalent:
\begin{enumerate}
\item $\underset{x\neq 0}{\underset{x\in \mathcal{X}}{\sup}}\frac{\norm{x}_\infty}{\norm{x}_\Phi}<\infty$.
\item $\underset{x\neq 0}{\underset{0\leq x\in \mathcal{X}}{\sup}}\frac{\norm{x}_\infty}{\norm{x}_\Phi}<\infty$.
\item $M$ is purely atomic, that is, $M=\underset{i=1}{\overset{N}\oplus} M_i$, $N\leq \aleph_0$ with $M_i$ is type $\mathrm{I}$ factor, and 
\begin{align*}
\tau(x)\geq \lambda (\underset{i=1}{\overset{N}\oplus} \Tr_i)(x), \quad 0\leq x\in M,
\end{align*}
holds true for some $\lambda>0$, where $\Tr_i$ is the canonical trace on $M_i$.
\end{enumerate}
\end{prop}

\begin{proof}
$(1)\Rightarrow(2)$ is trivial. For $(2)\Rightarrow(1)$, let $x\in\mathcal{X}$ and let $x=v\abs{x}$ be its polar decomposition. Since $\mathcal{X}$ is  a left ideal, $\abs{x}=v^*x\in\mathcal{X}$. Moreover, $\norm{x}_\infty=\norm{\abs{x}}_\infty$ and $\norm{x}_\Phi=\norm{\abs{x}}_\Phi$. Hence, the two suprema coincide.

\noindent $(3)\Rightarrow(1)$.
By $(4)\Rightarrow(5)$ of Theorem \ref{Thm:PurelyAtomic}, $\tau(e)\geq\lambda$ for every nonzero projection $e\in\mathcal{P}(M)$. Indeed, $\tau(e)\geq\lambda\Tr_i(ep_i)\geq\lambda$ for any $i$ with $ep_i\neq0$. Hence, Eq.  \eqref{Eq:BestConstant} with $\delta=\lambda$,
we have
\begin{align*}
\norm{x}_\infty\leq\Phi^{-1}(1/\lambda)\norm{x}_\Phi,\qquad x\in L^\Phi(M,\tau).
\end{align*}
In particular, this holds for all $x\in\mathcal{X}$, which gives $(1)$.

$(1)\Rightarrow(3)$. Let $C>0$ be such that $\norm{x}_\infty\leq C\norm{x}_\Phi$ for all $x\in\mathcal{X}$. We show that $L^\Phi(M,\tau)\subseteq M$. Let $x\in L^\Phi(M,\tau)$. Since $\mathcal{X}$ is dense in $L^\Phi(M,\tau)$, there exists a sequence $(x_n)$ in $\mathcal{X}$ with $\norm{x_n-x}_\Phi\rightarrow 0$. Note that $\norm{x_n-x_m}_\infty\leq C\norm{x_n-x_m}_\Phi$ for all $n,m$. Thus,$(x_n)$ is Cauchy in $M$ and hence, it converges in $\norm{\cdot}_\infty$ to some $y\in M$. As in the proof of $(6)\Rightarrow(7)$ of Theorem \ref{Thm:PurelyAtomic}, $x_n\rightarrow y$ and $x_n\rightarrow x$ in the measure topology, which is Hausdorff. Hence, $x=y\in M$. Thus, $L^\Phi(M,\tau)\subseteq M$, which is condition $(6)$ of Theorem \ref{Thm:PurelyAtomic}. Now condition $(3)$ follows from $(6)\Rightarrow(4)$ of Theorem \ref{Thm:PurelyAtomic}.
\end{proof}

\begin{rem}\label{Rem:XequalNp}
In conditions $(2)$-$(5)$ of Theorem
\ref{Thm:PurelyAtomicVonNeumannAlgebras} with $\Phi(t)=t^p$, the assumption that
$\mathcal{X}=\mathfrak{N}_p$ cannot be omitted to conclude
$\mathcal{X}$ admits an $M$-invariant algebraic complement in $L^p(M,\tau)$. Indeed, let $M=\mathbf{B}(\ell^2)$ and let $\tau$ be canonical trace on $\mathbf{B}(\ell^2)$ that takes value $1$ on each minimal projections. Clearly, $\tau(1)=\infty$. Moreover, $\mathfrak{N}_p=L^p(M,\tau)=\mathcal{S}^p(\ell^2)$, the Schatten $p$-operators. Note that the left ideal $\mathcal{X}$  of finite-rank operators is contained in $\mathfrak{N}_p$ and is $\norm{\cdot}_p$-dense in $L^p(M,\tau)$. Since $\norm{x}_\infty\leq\norm{x}_p$ for all $x\in \mathcal{X}$, the conditions of Proposition \ref{Prop:LeftIdealBound} hold (with $\Phi(t)=t^p$), whereas $\mathcal{X}\neq\mathfrak{N}_p$. Note that  $\mathcal{X}$ is not closed in $L^p(M,\tau)$ and by $(1)\Rightarrow (2)$ of Theorem \ref{Thm:PurelyAtomicVonNeumannAlgebras} admits no $M$-invariant algebraic complement in $L^p(M,\tau)$. Thus, none of the conditions of Theorem \ref{Thm:PurelyAtomicVonNeumannAlgebras} holds for this $\mathcal{X}$. In particular, $\mathcal{X}$ is not closed in $L^p(M,\tau)$ and admits no $M$-invariant algebraic complement in $L^p(M,\tau)$.
\end{rem}

We now have the following result, which characterizes finite-dimensional von Neumann algebras.  
\begin{thm}\label{Thm:PurelyAtomicVonNeumannAlgebras3}
Let $M$ be a finite von Neumann algebra with separable predual, equipped with a faithful normal tracial state $\tau$. Let $\Phi$ be an Orlicz function, and let $\mathcal{X}$ be a left ideal of $M$ contained in 
$\mathfrak{N}_\Phi:=M\cap L^\Phi(M,\tau)$ such that $\mathcal{X}$ is dense in $L^\Phi(M,\tau)$ with respect to $\norm{\cdot}_\Phi$.
Then the following statements are equivalent:
\begin{enumerate}
\item $\mathcal{X}$ admits an $M$-invariant algebraic complement in $L^\Phi(M,\tau)$, that is, there exists a linear subspace $\mathcal{Y}\subseteq L^\Phi(M,\tau)$ such that 
\begin{align*}
L^\Phi(M,\tau)=\mathcal{X}+\mathcal{Y}, \quad \mathcal{X}\cap\mathcal{Y}=\{0\}, \quad\text{and}\quad a\mathcal{Y}\subseteq\mathcal{Y}\,\text{ for all } a\in M.
\end{align*}

\item $\mathcal{X}=\mathfrak{N}_\Phi$ and $\mathfrak{N}_\Phi=L^\Phi(M,\tau)$.
			
\item 	$\mathcal{X}=\mathfrak{N}_\Phi$  and $\underset{x\neq 0}{\underset{0\leq x\in M}{\sup}}\frac{\norm{x}_\infty}{\norm{x}_\Phi}<\infty$.
			
\item 	$\mathcal{X}=\mathfrak{N}_\Phi$  and the reduced von Neumann algebra $eMe$ is finite-dimensional for every projection $e\in\mathcal{P}(M)$ with $\tau(e)<\infty$.
			
\item $\mathcal{X}=\mathfrak{N}_\Phi$  and $M$ is finite dimensional.
\end{enumerate}
\end{thm}

\begin{proof}
Since Theorem \ref{Thm:PurelyAtomicVonNeumannAlgebras} does not require any assumption on $\tau(1)$, it can be applied to $\tau$. Thus, Theorem \ref{Thm:PurelyAtomicVonNeumannAlgebras}  gives that $(1)\Leftrightarrow (2)\Leftrightarrow (3)\Leftrightarrow (4)$. Finally, $(4)$ with $e=1$ gives $(5)$, and $(5)$ trivially implies $(4)$. This completes the proof.
\end{proof}

\noindent\textbf{Acknowledgment}: The author acknowledges support from the INSPIRE Research Grant (Ref. No. DST/INSPIRE/04/2024/003789), DST, Government of India.

\noindent\textbf{Data Availability: }
Data sharing is not applicable to this article, as no new data
were created or analyzed in this study.	



\begin{thebibliography}{}


\bibitem[AEGV21]{AEGV21}
J. Alaminos, J.  Extremera, M. L. C. Godoy,  and A. R. Villena, 
{\emph{Hyperreflexivity of the space of module homomorphisms between non-commutative $L^p$-spaces}}, J. Math. Anal. Appl., {\bf 498}(2) (2021),  Paper No. 124964, 31 pp.
		
		
\bibitem[BDL90]{BDL90}  
{D. Buchholz, C. D'Antoni and R. Longo},
{\emph{Nuclear maps and modular structures. \rm{I}.\ {G}eneral properties}}, {J. Funct. Anal.}, \textbf{88}(2) (1990), {233--250}.
		


\bibitem[Da00]{Da00} 
H.~G. Dales, {\it Banach algebras and automatic continuity}, London Mathematical Society Monographs. New Series Oxford Science Publications, 24, Oxford Univ. Press, New York, 2000.



		
\bibitem[DM21]{DM21}
{D. De and K. Mukherjee}, 
{\emph{On symmetric embedding of type $\rm{I}$ factors}}, {J. Math. Phys.}, \textbf{62}(1) (2021), {Paper No. 013511, 28pp}.
		
		
		
\bibitem[DM22]{DM22}
{D. De and K. Mukherjee}, 
{\emph{A characterization of atomic von {N}eumann algebras}}, {Proc. Indian Acad. Sci. Math. Sci.}, \textbf{132}(2) (2022), {Paper No. 44, 9pp}.
		
		
\bibitem[DM23]{DM23}
{D. De and K. Mukherjee}, 
{\emph{On Bounded Coordinates in GNS spaces}},  {Dissertationes Math.}, \textbf{583} (2023), {1--106}.
		
		
		
		
		
\bibitem[FK86]{FK86}	
{T. Fack and H. Kosaki}, 
{\emph{Generalized {$s$}-numbers of {$\tau$}-measurable operators}},
		{Pacific J. Math.}, \textbf{123}(2) (1986), {269--300}.
		
		
		
		
		
		
		
		
		
		
\bibitem[FZ16]{FZ16}
{F. Fidaleo and L. Zsid\'o}, 
{\emph{Quantitative {BT}-theorem and automatic continuity for standard von {N}eumann algebras}}, {Adv. Math.}, \textbf{289} (2016), {1236--1260}.
		


\bibitem[GL20]{GL20}
S. Goldstein and L.~E. Labuschagne, {\it Notes on noncommutative $L^p$ and Orlicz spaces}, Wydawn. Uniw. \L\'odzkiego, \L\'od\'z, (2020).        
		
\bibitem[Ha75]{Ha75}
{U. Haagerup}, 
{\emph{The standard form of von {N}eumann algebras}}, {Math. Scand.}, \textbf{37}(2) (1975), {271--283}.
		
		
	
		
		
		
		
		
\bibitem[Hi21]{Hi21}
{F. Hiai}, 
{\emph{Lectures on selected topics in von {N}eumann algebras}},
		{EMS Series of Lectures in Mathematics}, {EMS Press, Berlin}, (2021).
		
		
		
		
		
\bibitem[HKS24]{HKS24}  	
{K. H. Han, S. H. Kye and E. St\o rmer},
        {\emph{Infinite dimensional analogues of {C}hoi matrices}}, {J. Funct. Anal.}, \textbf{287}(8) (2024), {Paper No. 110557, 23}.
	







\bibitem[Ka67]{Ka67}
R. V. Kadison, \emph{Problems on von Neumann algebras}., The Baton Rouge conference on operator algebras, Baton Rouge, (unpublished) (1967). 


\bibitem[KR97]{KR97}
R.~V. Kadison and J.~R. Ringrose, {\it Fundamentals of the theory of operator algebras. Vol. II}, corrected reprint of the 1986 original, 
Graduate Studies in Mathematics, 16, Amer. Math. Soc., Providence, RI, 1997.


\bibitem[Ro83]{Ro83}
J.~L. Romero, \emph{When is $L\sp{p}(\mu )$\ contained in $L\sp{q}(\mu )$?}, Amer. Math. Monthly {\bf 90}(3) (1983), 203--206.





        
		
\bibitem[SZ19]{SZ19}
{\c S.~V. Str\u atil\u a and L. Zsid\'o}, 
{\emph{Lectures on von Neumann Algebras}}, {Cambridge IISc Series},
		{Second edition}, {Cambridge University Press}, {(2019)}.
		


\bibitem[Vi85]{Vi85}
A. Villani, \emph{Another note on the inclusion $L^p(\mu)\subset L^q(\mu)$}, Amer. Math. Monthly {\bf 92}(7) (1985),  485--487.
		
\end{thebibliography}
\end{document}